\documentclass[a4paper]{amsart} 
\usepackage{amsmath}

\usepackage{amssymb}
\usepackage{verbatim}
\usepackage{amsthm}
\usepackage{mathrsfs}
\usepackage{graphicx}
\usepackage{mathtools}
\usepackage[colorlinks,pdfdisplaydoctitle,linkcolor=blue,citecolor=red]{hyperref}
\usepackage{caption}
\usepackage{subcaption}
\usepackage{url}
\usepackage{epstopdf}
\usepackage{pgf,tikz}
\usepackage{bbm}
\usepackage{extarrows}
\usetikzlibrary{arrows}

\renewcommand{\Delta}{\triangle}

\definecolor{darkblue}{rgb}{0,0,0.7}

\definecolor{darkgreen}{rgb}{0.01,0.75,0.24}

\def \Ee[#1]{\mathcal{E}^{\text{{#1}}}}
\def\R{\mathbf{R}}   

\def\bq{\mathbf{q}}

\def\pa[#1,#2]{\frac{\partial {#1}}{\partial {#2}} }

\def\idom[#1,#2,#3]{\int_{#1}\hspace{1pt} {#2} \hspace{1pt} \text{d}{#3}}
\def\res[#1,#2]{\left.{#1}\right|_{#2}}

\def\gt{\rightarrow}
\def\lgt{\downarrow}

\def\var[#1,#2]{\langle \delta \mathcal{E}^{\text{{#1}}}({#2}),v\rangle}
\def\vars[#1,#2,#3]{\langle \delta^2\mathcal{E}^{\text{{#1}}}({#2})v,{#3}\rangle}
\def\vard[#1,#2,#3,#4]{\langle \delta\mathcal{E}^{\text{{#1}}}({#2})-\delta\mathcal{E}^{\text{{#3}}}({#4}),v\rangle}

\def\F{\mathcal{F}}

\def\U{\mathcal{U}}

\newcommand{\bSig}{\boldsymbol{\Sigma}}

\newcommand{\bF}{\mathbf{F}}

\newcommand{\bG}{\mathbf{G}}
\newcommand{\bA}{\mathbf{A}}
\newcommand{\bB}{\mathbf{B}}

\newcommand{\A}{\mathcal{A}}

\newcommand{\eps}{\varepsilon}

\newcommand{\bI}{\mathbf{I}}

\newcommand{\bU}{\mathbf{U}}
\newcommand{\bu}{\mathbf{u}}
\newcommand{\bp}{\mathbf{p}}
\newcommand{\bY}{\mathbf{Y}}
\newcommand{\bQ}{\mathbf{Q}}

\newcommand{\be}{\begin{equation}}
\newcommand{\en}{\end{equation}}
\newcommand{\ben}{\begin{equation*}}
\newcommand{\enn}{\end{equation*}}
\newcommand{\bea}{\begin{aligned}}
\newcommand{\ena}{\end{aligned}}

\def\ba#1\ena{\begin{align}#1\end{align}}
\def\ban#1\enan{\begin{align*}#1\end{align*}}

\theoremstyle{plain}
\newtheorem{thm}{Theorem}[section]

\newtheorem{lem}[thm]{Lemma}
\newtheorem{cor}[thm]{Corollary}

\newtheorem{assumption}[thm]{Assumptions}
\newtheorem{assumptions}[thm]{Assumptions}

\theoremstyle{remark}
\newtheorem{rem}[thm]{Remark}

\numberwithin{equation}{section}

\begin{document}
\title[Bernstein-von Mises Theorem for nonlinear Bayesian inverse problems]{On the Bernstein-Von Mises Theorem for High Dimensional Nonlinear Bayesian Inverse Problems}

\begin{abstract}
We prove a Bernstein-von Mises theorem for a general class of high dimensional nonlinear Bayesian inverse problems in the vanishing noise limit. We propose a sufficient condition on the growth rate of the number of unknown parameters under which the posterior distribution is asymptotically normal. This growth condition is expressed explicitly in terms of the model dimension, the degree of ill-posedness of the inverse problem and the noise parameter. The theoretical results are applied to a Bayesian estimation of the medium parameter in an elliptic problem.
\end{abstract}

\author[Y. Lu]{Yulong Lu}
\address[Y. Lu]{Mathematics Institute, University of Warwick, Coventry CV4 7AL, UK}
\email{Yulong.Lu@warwick.ac.uk}

\date{\today}

\maketitle

\section{Introduction}
\subsection{Background and Aim}
In Bayesian statistics, the famous Bernstein-von Mises theorem \cite{VonMises31, le2012asymptotic} states that the posterior distribution in a Bayesian procedure is asymptically a normal distribution when the sample size tends to infinity or when the noise level tends to zero. To be concrete, consider the parametric case where one observes a sequence of independent and identically distributed random samples $Y_n$ of size $n$ from some distribution $P_f$, with $f$ belonging to some finite dimensional parameter space $\F$. Let $\Pi$ be a prior distribution on $f$ and denote by $\Pi(f | Y_n)$ the resulting posterior distribution. The Bernstein-von Mises theorem asserts that under the frequentist assumption that $Y_n$ is generated from some true parameter $f_0\in \F$, as $ n\gt \infty$, 
\be\label{eq:bvm-intro}
d_{TV} \Big(\Pi(f | Y_n), N(\Delta_n, \frac{1}{n} \mathcal{I}^{-1}_{f_0})\Big) \xrightarrow{P_{f_0}} 0,
\en 
where $\Delta_n$ is an efficient estimator for $f$, $\mathcal{I}_{f}$ is the Fisher information matrix of $P_f$ and $d_{TV}$ represents the total variation distance. The BvM theorem is important at least for two reasons. First, the BvM result \eqref{eq:bvm-intro} directly implies that the posterior contracts around the truth with the rate $\mathcal{O}(n^{-1/2})$. Second, Bayesian credible sets, which are sets that have fixed proportion of the total mass under the posterior measure, are asymptotically equivalent to frequentist confidence intervals, whereby the estimation of the latter can be realised by making use of the computational power of Markov chain Monte Carlo algorithms. 

In recent years, there has been growing interest in the areas of inverse problems and uncertainty quantification.  
 Inverse problems concern converting observed data into information about the quantities of interest which are not observed directly. Solving inverse problems is more difficult than solving the underlying forward problems because inverse problems are often ill-posed, meaning that they may have no (or multiple) solutions and that the inversion process is unstable. The latter scenario typically arises when the observations are contaminated by noises. Quantifying the uncertainty \cite{biegler2011large} associated with inevitable noises and inaccuracy of mathematical modelling is becoming an indispensable component in the resolution of inverse problems. The Bayesian approach \cite{kaipio2006statistical, DS15, S10a} provides a natural framework for doing so. It treats unknown quantities as random variables and updates the prior knowledge about unknowns by blending prior beliefs with data. This probabilistic approach offers several advantages over optimisation-based regularisation methods \cite{engl1996regularization}: On the one hand, the posterior distribution, the solution of Bayesian inverse problems, not only leads to a family of estimators, such as the posterior mean and the maximum posterior estimators, but also allow to quantify uncertainty inherent in Bayesian inferences; On the other hand, seeking solutions in wider probability spaces rather than in state spaces has a  ``stabilising" effect in the sense that posterior measure usually depends continuously on the data in some suitably sense.
  
 In this paper, we aim to prove the Berstein-von Mises theorem for Bayesian inferences of a general class of high dimensional nonlinear inverse problems in the limit of small noise. By ``high dimensionality", we mean that the number of the unknown parameters of the inverse problem is allowed to grow with the decreasing noise level. We want to develop some sufficient growth condition on the dimension of the model parameter under which one can observe the asymptotic normality of the posterior distribution.


\subsection{Relevant Literature and Contribution}

The reliability of a Bayesian approach can be assessed by the asymptotic performance of the posterior measure. This is formalised in the notion of {\em posterior consistency}. The posterior consistency refers to the contraction of the posterior distribution around the truth in the limit of large data size or small noise. 
In the context of Bayesian inverse problems, many posterior contraction results have been obtained for different problems, but most of them are focused on linear models; see e.g. \cite{KVZ11, KVZ14, ALS13, ASZ14} for linear Bayesian inverse problems with Gaussian priors and \cite{R13} for non-conjugate priors. Recently, Vollmer \cite{V13} proved a posterior contraction result for a specific nonlinear Bayesian inverse problem arising from subsurface flow modelling, based on some stability estimates of the forward problem and nonparametric regression results.   

The Bernstein-von Mises theorem is a classical result in parametric statistical models (cf. \cite{V00, le2012asymptotic}). The last decade has seen a significant development of BvM theorems in infinite dimensional statistical models. For nonparametric statistical models, whether the BvM theorem hold depends crucially on the topology under which the convergence of probability measures is taken. In fact, it was observed in \cite{C93, F99} that the nonparametric BvM theorem can not be true in the strictly $L^2$-sense. Nevertheless, recent work by Castillo and Nickl \cite{castillo2013nonparametric, castillo2014bernstein} showed that the BvM theorem does hold when the metric of convergence is relaxed to be some suitable weak convergence (essentially convergence in 1-Wasserstein distance) of probability measures in a larger space than $L^2$ (in fact some negative Sobolev space). Also, semiparametric BvM results had been obtained by many authors in various types of statistical models, see e.g. \cite{shen2002asymptotic, kim2006bernstein,rivoirard2012bernstein,bickel2012semiparametric,castillo2012semiparametric, castillo2015bernstein,PS14}. 

The type of BvM results of our particular interest is the case where the dimension of parameters in the underlying statistical model can grow with the increasing sample size or decreasing noise level. Ghosal \cite{Ghosal99, Ghosal00} studied the behaviour of the posterior distribution when the number of parameters tend to infinity for a linear regression model and exponential families. Under certain growth condition on the model dimension, he proved that the posterior is approximately a normal distribution. Similar results have also been obtained in other statistical models, see \cite{bontemps2011bernstein} for Gaussian regression with increasing number of regressors,  \cite{boucheron2009bernstein} for discrete probability distributions, \cite{clarke2010reference} for exponential families with reference priors. It is worth to mention that Belloni and Chernozhukov \cite{BC14} revisited exponential families and improved Ghosal's growth rate condition by a logarithm factor. Recent results by Spokoiny and Panov \cite{PS14, spokoiny2013bernstein}  show that to a certain extent this new growth condition seems to be sharp. 

Despite the extensive study of the BvM theorem in the statistical community, there are only few BvM results in Bayesian inverse problems. Knapik et al \cite{KVZ11} proved a BvM theorem for linear projections of the posterior distribution of Bayesian linear problems with Gaussian priors. In particular, they showed that the Bayesian credible set asymptotically agrees with the frequentist confidence region in the undersmoothing case where the prior is less regular than the truth,  whereas oversmoothing leads to zero frequentist coverage. For nonlinear Bayesian inverse problems, Lu et al \cite{lu2016gaussian} recently proved a parametric BvM result, which characterises the Gaussian approximations of posterior measure with respect to the Kullback Leibler divergence. 

To the best of our knowledge, the BvM theorem for nonlinear Bayesian inverse problems with infinite number of unknowns has not been addressed in the literature. The present work gives a first attempt in this direction. Our aim is to generalise Ghosal's previous results \cite{Ghosal99, Ghosal00} on exponential families and linear regressions to a Bayesian inverse problem setting. In particular, we focus on a general finite dimensional Bayesian inverse problem and study the behaviour of the posterior distribution in the limit of small noise where the number of unknowns is allowed to grow to infinity.  
 The finite dimensional forward problem of consideration is strongly motivated by the discretisation of partial differential equation (PDE) models; the dimension of the problem is then associated to the number of the mesh points or basis functions of a discretisation scheme. The solution of the forward problem approximates to the solution of the original PDE as the model dimension grows up to infinity. The ill-posedness of the corresponding inverse problem, the focus of our attention, is quantified in terms of the decaying rate of the singular values of the linearisation of the nonlinear forward model. By prescribing an appropriate prior on the unknown parameter, we show that, under certain growth constraint on the dimension, the posterior distribution can be approximated by a normal distribution centring around the truth in the total variation distance (see Theorem \ref{thm:main}). The growth condition depends on the dimension, the decaying rate of singular values and the decreasing noise level. This result also improves our previous BvM result \cite{lu2016gaussian} for Bayesian inverse problems with fixed finite dimension. Our proof follows closely Ghosal's idea \cite{Ghosal99}, but we simplify his arguments substantially by taking advantage of the specific structure of the posterior in the inverse problem setting. We also verify the assumptions that have been made leading to asymptotic normality in an inverse medium problem.
 
 We emphasise that there is still an essential gap between our BvM result and the nonparametric BvM theorem in \cite{KVZ11} since we are unable to justify the frequentist coverage of credible sets defined on the full (infinite dimensional) parameter space. It is expected that the asymptotic normality of the posterior is still valid under certain regular linear projections \cite{Nickl}. Due to the nonlinearity and the ill-posedness of the inverse problem, justifying this is highly non-trivial and is far beyond the scope of the present paper. Finally, it is also unclear whether the weak BvM theorems obtained in \cite{castillo2013nonparametric, castillo2014bernstein} still hold in the context of nonlinear Bayesian inverse problem. Nevertheless, our results may give a hint about how to choose the weak topology that leads to BvM theorems. This will be examined carefully in future work.
 
\subsection{Structure}
This paper is organised as follows. In Section \ref{sec:notation}, we define some useful notations that will be used throughout the paper. In Section \ref{sec:setup} we make some assumptions on the forward problem as well as the prior and then state the main BvM results in Theorem \ref{thm:main} and Corollary \ref{cor:bvm-Pi}. Proofs of the main results are provided in Section \ref{sec:proofmain}. Section \ref{sec:appl} discusses an application of previous theoretical results in an inverse medium problem. 

\subsection{Notations}\label{sec:notation}
We denote by $|f|$ the standard Euclidean norm of a vector $f\in \R^d$ and by $\|\bA\|$ the operator norm of a matrix $\bA$, i.e. $\|\bA\| = \sup_{|f|\leq 1} |\bA f|$. Let $\langle f, g \rangle$ be the inner product of vectors $f, g \in \R^d$. We denote by $\sigma_{\min}(\bA)$ and $\sigma_{\max}(\bA)$ (or $\lambda_{\min}(\bA)$ and $\lambda_{\max}(\bA)$) the minimum and maximum singular values (or eigenvalues) of $\bA$ respectively.  Let $\text{Vol}(d)$ denote the volume of the $d$-dimensional unit ball.
 We use $C$ to denote a generic constant independent of $d$ and $n$ and we write $a \lesssim b$ when $a \leq C b$. For a sequence of random variables $\{X_n\}$, we write $X_n = \mathcal{O}_p(\delta_n)$ (or $X_n = o_p(\delta_n)$) to indicate that the sequence $X_n/\delta_n $ is bounded in probability (or $X_n/\delta_n \gt 0$ in probability). 

\section{Set-up}\label{sec:setup}
Consider the Bayesian inverse problem of estimating $f\in \R^d$ from noisy measurement $Y_n \in\R^{d}$ with 
\be\label{eq:ip}
Y_n = G(f) + \frac{1}{\sqrt{n}}\eta.
\en
Here the operator $G: \R^d \gt \R^d$ is the forward map which is considered to be nonlinear. We also assume that $\eta$ is a standard normal distributed random variable in $\R^d$. The forward problem \eqref{eq:ip} in our mind comes from a finite dimensional approximation to some infinite dimensional problem. Typically, in the case where the underlying forward model is governed by a PDE, the forward operator $G$ in \eqref{eq:ip} may be obtained from some finite dimensional discretisation of the PDE, for example, through finite difference methods or Galerkin methods. The data $Y_n$ could be generated through noisy observations of the solution to the PDE at $d$ spacial/temporal positions. The standard normal assumption of the noise $\eta$ is of particular interest to us because it can be viewed as a discrete analogue of a white noise process. 

Given a prior to be defined below, we are interested in the asymptotic performance of the resulting posterior as $n\gt \infty$ under the frequentist assumption that $Y_n$ is generated from some truth $f_0$.  We assume that $f_0$ lies in a compact subset $\mathcal{F}$ on $\R^d$. We also make the following assumptions on the forward operator $G$. 
\begin{assumption}\label{as:G}

\item[(A1)] For every $f\in \F$, $G(f)$ is differentiable and its derivative matrix $\nabla G(x)$ is invertible. Moreover, there exists positive constants $\sigma_0$ and $\sigma(d)$ such that
$$
\sigma(d)^{-1} \leq \sigma_{\min}(\nabla G(f))\leq  \sigma_{\max}(\nabla G(f))  \leq\sigma_0^{-1} \quad \text{ for all } f\in \F.
$$

\item[(A2)] There exist constants $C,\tilde{C} > 0$ such that 
$$
\tilde{C}|G(f_1) - G(f_2)| \leq |f_1- f_2| \leq C \sigma(d)|G(f_1) - G(f_2)|  \quad \text{ for all } f_1, f_2\in \F.
$$

\item[(A3)] There exists a constant $C > 0$ such that 
$$
|G(f_1) - G(f_2) - \nabla G(f_2) (f_1  - f_2)| \leq C |f_1 - f_2|^2 \quad \text{ for all } f_1, f_2\in \F.
$$
\end{assumption}

\begin{rem} 
\item [(i)] The assumption (A1) is equivalent to that
$$
\sigma_0 \leq \sigma_{\min}((\nabla G(f))^{-1})\leq  \sigma_{\max}((\nabla G(f))^{-1})  \leq \sigma(d)\quad \text{ for all } f\in \F.
$$
In particular, it also implies that $\|\nabla G(f)\|\leq \sigma_0^{-1}$.  The dimension-increasing constant $\sigma(d)$ bounds the growth of singular values of $\nabla G^{-1}$,   characterising the degree of ill-posedness of the infinite dimensional inverse problem. For mild ill-posed inverse problems $\sigma(d) \sim \mathcal{O}(d^s)$ with some $s> 0$ when $d \gg 1$, and for severely ill-posed inverse problems $\sigma(d) \sim \mathcal{O}(\exp(d^s))$.

\item[(ii)] Assumptions (A2) and (A3) are both regularity assumptions on the forward operator. In particular, (A3) holds if $G$ has bounded second derivatives. However, we avoid imposing constraints on the second order derivatives of $G$ since they are hard to evaluate or estimate in practice. It is also worth to remark that the constants $C, \tilde{C}$ appearing in the above assumptions are independent of a concrete $f$, but they may depend on the choice of $\F$. 

\item[(iii)] Assumptions (A3) and the growth condition on the singular values in (A1) imply that the estimate of (A1) holds when $|f_1 - f_2|$ is small. However, we assume that this is valid for all $f_1, f_2 \in \F$. This assumption is not strong as it looks since we will show that it is satisfied in many applications.  
\end{rem}

 We denote by $\Pi(d f)$ the prior distribution which has a Lebesgue density $\pi(f)$ on $\R^d$. Moreover, we assume $\Pi$ satisfies the following. 

\begin{assumptions}\label{as:prior}
\item[(i)] $\Pi$ is supported on the set $\F$. 

\item[(ii)] There exists $C>0$ such that 
$$
\sup_{f\in \F}\log (\pi(f)/ \pi(f_0)) \leq Cd.
$$

\item[(iii)] For any $\delta > 0$, there exists $C_\delta > 0$ such that 
$$
|\log \pi(f) - \log \pi(f_0)| \leq C_\delta \sqrt{d} |f - f_0|
$$
when $| f - f_0| \leq \delta\sqrt{d}$.
\end{assumptions}

\begin{rem}\label{rem:prior}
Assumption \eqref{as:prior} (ii) means that the value of prior density at the truth is not exponentially smaller than its value at any point in $\F$. Such assumption was also used in \cite{BC14}.  Assumption \eqref{as:prior} (iii) is a continuity condition for the log prior density near the truth. Similar assumptions was required in \cite{Ghosal99, Ghosal00, BC14}. In particular, these assumptions are fulfilled when the prior $\Pi$ on $f$ is a product of independent priors on each component $f_i$, i.e. $\pi(f)  = \Pi_{i=1}^d \pi(f_i)$, and each $\pi_i$ satisfies 
$
|\log \pi_i(f_i) - \log \pi_i(f_{0,i})| \leq C $ for any $f_i$ and 
$ |\log \pi_i(f_i) - \log \pi_i(f_{0,i})| \leq C_\delta |f_i - f_{0,i}| 
$
when $|f_i - f_{0,i}| \leq \delta$ for some $\delta  > 0$.
\end{rem}

Given the prior $\Pi$, the posterior distribution given the observed data $Y_n$, is defined as
$$
\Pi_n (d f) = \frac{\exp\Big(-\frac{n}{2} |Y_n - G(f)|^2\Big)}{\int_{\F} \exp\Big(-\frac{n}{2} |Y_n- G(\tilde{f})|^2\Big) \Pi(d\tilde{f})}  \Pi(d f).
$$
By dividing the common factor $\frac{n}{2}|Y_n|^2$ in the exponential, we can write the density of the posterior $\Pi_n$ as 
$$
\pi_n(f) = \frac{\exp\Big(n\big( \langle Y_n, G(f) \rangle - \frac{1}{2} |G(f)|^2\big)\Big)}{\int_{\F}  \exp\Big(n\big( \langle Y_n, G(\tilde{f}) \rangle - \frac{1}{2} |G(\tilde{f})|^2\big)\Big) \pi(\tilde{f})d \tilde{f}} \pi(f).
$$

We aim to prove the Bernstein-von Mises theorem for the posterior measure $\Pi_n$ when the dimension $d$ and the noise parameter $n$ increase to infinity simultaneously.  To this end,  it is more convenient to analyse the posterior distribution in the local parameter space around $f_0$. More specifically, let us define $\U := \sqrt{n} (\F - f_0)$. Let $u :=  \sqrt{n}(f - f_0) \in \U$. Then the posterior distribution of $u$, denoted by $\Pi_n^\ast$, has the density
$$
\pi_n^\ast (u) = \frac{\exp\Big(n\big( \langle Y_n, G(f_0 + n^{-\frac{1}{2}} u) \rangle -  \frac{1}{2}|G(f_0 + n^{-\frac{1}{2}} u)|^2\big)\Big)  \pi(f_0 + n^{-\frac{1}{2}} u)}{\int_{\U}  \exp\Big(n\big( \langle Y_n, G(f_0 + n^{-\frac{1}{2}} \tilde{u}) \rangle -  \frac{1}{2} |G(f_0 + n^{-\frac{1}{2}} \tilde{u})|^2\big)\Big) \pi(f_0 + n^{-\frac{1}{2}} \tilde{u})d \tilde{u}}.
$$

 Formally, the asymptotic normality of the posterior can be read from the linear expansion of the log likelihood ratio. 
In fact, for any $u\in \U$, let us define the log likelihood ratio
$$
L_n(u) := n\big( \langle Y_n, G(f_0 + n^{-\frac{1}{2}} u) - G(f_0) \rangle - \frac{1}{2} ( |G(f_0 + n^{-\frac{1}{2}} u) |^2-  |G(f_0)|^2)\big)
$$
and the shifted likelihood function
$
Z_n (u) = \exp(L_n(u)).
$
Then we can rewrite the posterior density $\pi_n^\ast$ in terms of $Z_n$, namely,
\be\label{eq:post1}
\pi^\ast_n (u) = \frac{Z_n(u) \pi(f_0 + n^{-\frac{1}{2}} u)}{\int_{\U}Z_n(\tilde{u})\pi(f_0 + n^{-\frac{1}{2}} \tilde{u}) d\tilde{u}}.
\en
Let us set $\bSig := (\nabla G(f_0)^{T} \nabla G(f_0))^{-1}$ and $\Delta_n := \bSig \nabla G(f_0)^T \eta$. Then by the normal assumption on the noise $\eta$, we have that $\Delta_n \sim N(0, \bSig)$. 
We also define $\tilde{Z}_n(u) = \exp(\tilde{L}_n(u))$ with the exponent 
$$
\bea
\tilde{L}_n(u) & := 2\langle u , \bSig^{-1} \Delta_n\rangle - |\bSig^{-\frac{1}{2}} u|^2\\
& = 2\langle \eta, \nabla G(f_0) u \rangle - |\bSig^{-\frac{1}{2}} u|^2.
\ena
$$
Recall that $\phi(u; m, \bSig)$ denotes the probability density function of the normal distribution $N(m, \bSig)$. 
Then it is easy to see that $\phi(u;\Delta_n, \bSig) = \tilde{Z}_n(u)/ \int_{\R^d}  \tilde{Z}_n(\tilde{u})d\tilde{u}$.

After completing square the log likelihood ratio $L_n$ can be expressed as
$$
\bea
L_n(u) & = n\big( \langle \bar{Y} - G(f_0), G(f_0 + n^{-\frac{1}{2}} u) - G(f_0) \rangle - \frac{1}{2} |G(f_0 + n^{-\frac{1}{2}} u) - G(f_0)|^2\big)\\
 & =  \langle \eta, \sqrt{n}(G(f_0 + n^{-\frac{1}{2}} u) - G(f_0) )  \rangle - \frac{n}{2}|G(f_0 + n^{-\frac{1}{2}} u) - G(f_0)|^2.
\ena
$$
By expanding the function $G(f_0 + n^{-\frac{1}{2}} u)$ around the origin up to the first order, one can see that at least locally
$$
L_n(u) \approx \langle \eta, \nabla G(f_0) u\rangle -\frac{1}{2} |\nabla G(f_0) u|^2 = \tilde{L}_n(u).
$$
This combing with $\pi(f_0 + n^{-\frac{1}{2}} u) \approx \pi(f_0)$ implies that $\pi^\ast_n(u) \approx \phi(u; \Delta_n, \bSig)$. The formal calculations above can be made rigorous under certain growth condition on the dimensionality $d$ with respect to $n$.  

Given a fixed $K> 0$, let us define
\be\label{eq:Kd}
K(d) := K\sigma(d) \sqrt{d (\log(d) + \log \sigma(d))}
\en
 and let $\delta_n :=\sqrt{\frac{d}{n}} K^2(d)$. Our main result is the following Bernstein-von Mises theorem.

\begin{thm}\label{thm:main}
Let Assumptions \ref{as:G} and Assumptions \eqref{as:prior} be satisfied. If $\delta_n \gt 0$ as $n \gt \infty$.  Then we have
$$
d_{TV}\Big(\Pi_n^\ast, N(\Delta_n, \bSig) \Big) = \int_{\U} |\pi^\ast_n(u) - \phi(u; \Delta_n, \bSig)| du = \mathcal{O}_p(\delta_n)
$$
as $n \gt \infty$. 
\end{thm}

Let $\overline{\Delta}_n := n^{-\frac{1}{2}} \Delta_n + f_0$. Then $\overline{\Delta}_n \sim N(f_0, n^{-\frac{1}{2}} \bSig)$.  Since the total variation norm is invariant under the bijection $f\mapsto \sqrt{n}(f - f_0)$, we state the following BvM result for the original posterior measure $\Pi_n$, as a corollary of Theorem \ref{thm:main}. 
\begin{cor}\label{cor:bvm-Pi}
Under the same assumption as in Theorem \ref{thm:main}, we have
$$
d_{TV}\Big(\Pi_n, N(\overline{\Delta}_n, n^{-1}\bSig) \Big) = \int_{\F} |\pi_n(f) - \phi(f; \overline{\Delta}_n, n^{-1}\bSig)|df = \mathcal{O}_p(\delta_n)
$$
as $n\gt \infty$.
\end{cor}

\begin{rem}
As an important consequence of Corollary \ref{cor:bvm-Pi}, the Bayesian credible set is asymptotically identical to the standard frequent conference interval. To be more precise, given a credible level $\alpha \in (0,1)$, let $C_{n, \alpha}$ be the credible set such that $\Pi_n (C_{n ,\alpha}) = 1 - \alpha$. If $\delta_n \gt 0$ as $n \gt \infty$, then it follows from Corollary \ref{cor:bvm-Pi} and the fact that $\Delta_n = \mathcal{O}_p(\sqrt{d\|\bSig\|}) = o(n^{\frac{1}{2}})$ that 
$$
P_{f_0}^n (f_0 \in C_{n, \alpha}) \gt 1 - \alpha.
$$
\end{rem}

\begin{rem}It is worth to make a comment on the growth requirement on the dimension $d$, that is,  $\delta_n =\sqrt{d/n} K^2(d) \gt 0$. In fact, this condition essentially agrees with the conditions proposed by Ghosal \cite{Ghosal99, Ghosal00} for proving the BvM results for linear regression models (see Assumption (A4) in \cite{Ghosal99}) and exponential families (see Condition (R) in \cite{Ghosal00}). The only difference is that his conditions were expressed in terms of the norm of the Fisher information matrix and our condition is written in terms of $\sigma(d)$. Recently, Belloni and Chernozhukov \cite{BC14} relaxed Ghosal's growth requirement for exponential families by removing the logarithm factors and obtained that the BvM theorem holds when $d^3/n \gt 0$.  It has been shown in \cite{PS14} that this growth is indeed sharp for a specific i.i.d smooth statistical model. This suggests that the sharp growth condition for the BvM for Bayesian inverse problems might be $\sqrt{\frac{d^3}{n}}\sigma(d)^2\gt 0$, however we are unable to justify this yet. This is to be investigated in future work.

\end{rem}

\section{Proof of the Main Result} \label{sec:proofmain}

\subsection{Lemmas} 
The proof of Theorem \ref{thm:main} follows directly from a series of preliminary lemmas as we establish now.  
The first lemma in the following gives an estimate for the tail probability of the normal distribution $N(\Delta_n, \bSig)$.

\begin{lem}\label{lem:normaltail}
Given any $c > 0$, there exists $K> 0$ such that, with probability tending to one, 
$$
\int_{|u| \geq K \sigma(d)\sqrt{d}} \phi(u; \Delta_n, \bSig) du \leq e^{-cd}.
$$
\begin{proof}
First notice that $\Delta_n \sim N(0, \bSig)$. Then by the definition of $\bSig$ and Assumption \ref{as:G} (A1), we can obtain from Chebyshev's inequality that $|\Delta_n| = \mathcal{O}_p(\sqrt{|\bSig| d}) = \mathcal{O}_p(\sigma(d)\sqrt{d}) $. Therefore with probability tending to one, the following holds when $K$ is sufficiently large:
$$
\bea
& \int_{|u| \geq K \sigma(d)\sqrt{d}} \phi(u; \Delta_n, \bSig) du \leq \int_{|u| \geq K \sigma(d)\sqrt{d}} \phi(u; 0, \bSig) du \\
& = \int_{|u|\geq K\sqrt{d}} \phi(u; 0, \bI_d) du \\
& \leq e^{-cd}.
\ena
$$
Note that the last inequality follows from the classical tail inequality of standard normal distribution.
\end{proof}
\end{lem}

The next lemma states that the tail probability of the posterior distribution outside some large ball is negligible. 

\begin{lem}\label{lem:pitail}
 Let $B_K$ be a centred ball of radius $K(d) = K\sigma(d) \sqrt{d (\log(d) + \log \sigma(d))}$. Let $|\eta| \leq \tilde{C}\sqrt{d}$. Then there exists $K$ sufficiently large such that when $n$ is large, 
$$
 \int_{B_K^c \cap \U} \pi^\ast_n(u) du \leq \exp\Big(- d \big(\log(d) + \log \sigma(d)\big)\Big).
 $$
\begin{proof}
From the definition of $\pi^\ast_n$ in \eqref{eq:post1} and Assumption \ref{as:prior} (ii), it suffices to show that 
$$
\Big(\int_{\U}Z_n(\tilde{u}) d\tilde{u} \Big)^{-1} \int_{B_K^c \cap \U} Z_n(u)  du  \leq \exp\Big(- d \big(\log(d) + \log \sigma(d)\big)\Big).
$$
For doing so, we first prove an upper bound for the integral of $Z_n$ over $B_K^c \cap \U$. Recall that $Z_n(u) = \exp(L_n(u))$ and that 
$$
L_n(u) = \langle \eta, T_n(u)  \rangle - \frac{1}{2}|T_n(u)|^2
$$
where $T_n(u) = \sqrt{n} (G(f_0 + n^{-\frac{1}{2}} u) - G(f_0)) $. According to Assumption \ref{as:G} (A1), $|T_n(u)| \geq  C\sigma(d)^{-1}|u| \geq C\sigma(d)^{-1} K(d)$ when $u\in B_{K}^c \cap \U$. From the assumption that $\eta \leq \tilde{C}(\sqrt{d})$ and the definition of $K(d)$, we can choose $K$ to be sufficiently large so that $|\eta | \leq \frac{1}{2} |T_n(u)|$. With such $K$ being fixed, we have that
\be\label{eq:upperZ}
\begin{aligned}
 \int_{B_K^c \cap \U} Z_n(u)  du & \leq  \int_{B_K^c \cap \U} \exp \Big(\langle \eta, T_n(u)  \rangle  - \frac{1}{2}|T_n(u)|^2\Big) du \\
 & \leq  \int_{B_K^c \cap \U} \exp \Big(- \frac{1}{4}|T_n(u)|^2\Big) du\\
 &\leq \int_{B_K^c \cap \U} \exp \Big(- \frac{1}{4C \sigma(d)^2}|u|^2\Big) du\\
&  \leq C\exp\Big(-\frac{K^2(d)}{\sigma^2(d)}\Big)\times  (\sigma(d))^d\\
& = C\exp\Big(-\frac{K^2(d)}{\sigma^2(d)} + d\log \sigma(d)\Big).
 \end{aligned}
\en
Next, we seek an lower bound for the integral $\int_{\U} Z_n(\tilde{u})d \tilde{u}$. In fact,
$$
\begin{aligned}
 \int_{\U} Z_n(\tilde{u})d \tilde{u} & \geq \int_{\{\tilde{u} \in \U: |\tilde{u}| \leq 1\}} \exp\Big(\langle \eta, T(\tilde{u} )  \rangle - \frac{1}{2}|T(\tilde{u} )|^2 \Big) d \tilde{u}\\
& \geq \int_{\{\tilde{u} \in \U: |\tilde{u}| \leq 1\}} \exp\Big(-|\eta| |T(\tilde{u})| -\frac{1}{2}|T(\tilde{u} )|^2  \Big)d \tilde{u}.
\end{aligned}
$$
Assumption \ref{as:G} (A2) yields that $|T(\tilde{u})| \leq C$ when $|\tilde{u}|\leq 1$. Moreover, from the definition of $\U$, the unit ball in contained in $\U$ when $n$ is sufficiently large.   As a consequence, when $n$ is large, 
\be\label{eq:lowZ}
\begin{aligned}
 \int_{\U} Z_n(\tilde{u})d \tilde{u}  & \geq \exp\Big(-\tilde{C}\sqrt{d} C - \frac{1}{2}C^2\Big) \text{Vol} (\{\tilde{u} \in \U: |\tilde{u}| \leq 1\})\\
& \geq \exp\Big(-\tilde{C}\sqrt{d} C - \frac{1}{2}C^2\Big) \text{Vol}(d)
 \end{aligned}
\en
The volume of a $d$-dimensional unit ball is 
$$
\text{Vol}(d) = \frac{\pi^{d/2}}{\Gamma(d/2 +1)}
$$
and from Stirling's formula for the Gamma function we know that when $d \gg 1$
$$
\text{Vol}(d) \sim \Big(\frac{2\pi e}{d}\Big)^{d/2}. 
$$
Combining this with \eqref{eq:lowZ} and relabelling the constants $C$ and $\tilde{C}$ yields 
$$
\int_{\U} Z_n(\tilde{u}) d\tilde{u} \geq C \exp(-\tilde{C}(\sqrt{d} + d\log{d})  ).
$$
 Finally, the lemma follows from the above lower bound and the upper bound \eqref{eq:upperZ}. 
\end{proof}
\end{lem}

The following posterior contraction result is a direct consequence of Lemma \ref{lem:pitail}.
\begin{cor}
Let $K(d)$ be defined as in Lemma \ref{lem:pitail}. Assume that $\eps_n:=n^{-\frac{1}{2}} K(d) \gt 0$ as $n\gt \infty$. Then with probability tending to one, the original posterior measure of $f$ contracts around the truth $f_0$ with rate $\eps_n$, i.e.  for any $M_n \gt \infty$, 
$$
\Pi_n \Big(\{f: |f - f_0| \geq M_n\eps_n\}\Big) \gt 0
$$
as $n \gt \infty$. 
\end{cor}

\begin{lem}\label{lem:diffZ}
 Let $B_K$ be a centred ball of radius $K(d)$. Let $|\eta| \leq \tilde{C}\sqrt{d}$. If $\delta_n = \sqrt{\frac{d}{n}}K^2(d) \gt 0$ when $n \gt \infty$, then 
\be\label{eq:diffZ}
\Big(\int_{\U} \tilde{Z}_n(u)du\Big)^{-1}\int_{B_K \cap \U}  \Big|Z_n(u) - \tilde{Z}_n(u)\Big| = \mathcal{O}(\delta_n).
 \en 
 \begin{proof}
 From the definition of $Z_n$ and $\tilde{Z}_n$, and by Assumption \ref{as:G} (A3), we have for $u\in B_K \cap\U$, 
 \be
 \begin{aligned}
 \big| \log Z_n(u) -  \log \tilde{Z}_n(u)\big| & = \Big| \langle \eta, \sqrt{n}\big(G(f_0 + n^{-\frac{1}{2}} u) - G(f_0) - n^{-\frac{1}{2}} \nabla G(f_0)u\big)\rangle \\
& \quad -\Big(n \big|G(f_0 + n^{-\frac{1}{2}} u) - G(f_0) \big|^2 - |\nabla G(f_0) u|^2\Big) \Big|\\
&\leq C\Big(n^{-\frac{1}{2}} |\eta| |u|^2 + n^{-1}|u|^4\Big)\\
& \leq C\Big(n^{-\frac{1}{2}} |\eta| K^2(d) + n^{-1}K^4(d)\Big)\\
& \leq C\delta_n. 
\end{aligned}
 \en
It follows that 
 $
Z_n(u) \leq \tilde{Z}_n(u) + e^{\delta_n}.
 $
 An application of the elementary inequality $|e^a - e^b| \leq |a - b| \max(e^a, e^b)$ gives  
 $$
 \bea
 \Big(\int_{\U} \tilde{Z}_n(u)du\Big)^{-1} \int_{B_K \cap \U}  \Big|Z_n(u) - \tilde{Z}_n(u)\Big| & \leq C \delta_n e^{\delta_n} \Big(\int_{\U} \tilde{Z}_n(u)du\Big)^{-1}\int_{B_K \cap \U} \tilde{Z}_n(u)du \\
& \leq C\delta_n e^{\delta_n}.
\ena
 $$
 This proves \eqref{eq:diffZ}.

 \end{proof}
\end{lem}

Finally, we recall the following useful lemma from \cite{Ghosal99}.
\begin{lem}\label{lem:bdf}
Let $f, g$ be two non-negative integrable function not identically zero on a measurable space $S$ and let $F\subset S$. Then 
$$
\int_{F} \Big|\frac{f}{\int f} - \frac{g}{\int g}\Big| \leq \frac{\int_{F^c} f}{\int f} +  \frac{\int_{F^c} g}{\int g} + 3\Big(\int g\Big)^{-1} \int_F \Big| f - g\Big|.
$$
\end{lem}

\subsection{Proof of Theorem \ref{thm:main}}\label{sec:proofmain}

\begin{proof}
Recall the definition of the ball $B_K = \{u: |u| \leq K(d)\}$ with $K(d)$ defined as \eqref{eq:Kd}. Then by the triangle inequality,
$$\bea
& \int_{\U} |\pi^\ast_n(u) - \phi(u; \Delta_n, \bSig)| du \leq \int_{B_K \cap \U}\Big| \frac{Z_n(u) \pi(f_0 + n^{-\frac{1}{2}} u) }{\int_{\U}Z_n(u) \pi(f_0 + n^{-\frac{1}{2}} u) du} - \frac{\tilde{Z}_n(u) \pi(f_0) }{\int_{\R^d} \tilde{Z}_n(\tilde{u})d\tilde{u}} \Big| d u\\ 
& + \int_{B_K^c \cap \U} \pi^\ast_n(u) d u + \int_{B_K^c \cap \U} \phi(u; \Delta_n, \bSig) d u
\ena
$$
According to Lemma \ref{lem:bdf}, the first term on the right side of above equation can be bounded from above by
$$
\bea
&  \int_{B_K^c \cap \U} \pi^\ast_n(u) d u + \int_{B_K^c \cap \U} \phi(u; \Delta_n, \bSig) d u\\
& + 3 \Big(\int_{\R^d} \tilde{Z}_n(\tilde{u})d\tilde{u}\Big)^{-1} \int_{B_K \cap \U}  \Big|Z_n(u) \pi(f_0 + n^{-\frac{1}{2}} u) - \tilde{Z}_n(u) \pi(f_0)\Big| du.
 \ena
$$
 In addition, thanks to the fact that $\eta = \mathcal{O}_p(\sqrt{d})$, the first two terms above can be made sufficiently small with high probability by Lemma \ref{lem:pitail} and Lemma \ref{lem:normaltail} respectively. Hence the theorem is proved if we can show that the last term above is small with high probability when $n\gt \infty$. In fact, the last term can be bounded from above as
\be\label{eq:ratio0}
\begin{aligned}
& \Big(\int_{\R^d} \tilde{Z}_n(\tilde{u})d\tilde{u}\Big)^{-1} \int_{B_K \cap \U}  \Big|Z_n(u) \pi(f_0 + n^{-\frac{1}{2}} u) - \tilde{Z}_n(u) \pi(f_0)\Big| du \\
& \leq \sup_{u\in B_K\cap \U}\Big|\frac{\pi(f_0 + n^{-\frac{1}{2}}u)}{\pi(f_0)} - 1\Big| \Big(\int_{\U} \tilde{Z}_n(u)du\Big)^{-1} \int_{B_K\cap \U} Z_n(u)du\\
& \quad\quad + \Big(\int_{\U} \tilde{Z}_n(u)du\Big)^{-1}\int_{B_K \cap \U}  \Big|Z_n(u) - \tilde{Z}_n(u)\Big| du .
\end{aligned}
\en
By Assumption \ref{as:prior} (iii), when $n$ is sufficiently large, 
\be\label{eq:ratio1}
\begin{aligned}
\sup_{u\in B_K\cap \U}\Big|\frac{\pi(f_0 + n^{-\frac{1}{2}}u)}{\pi(f_0)} - 1\Big| & \leq \sup_{u\in B_K\cap \U} 2 |\log(\pi(f_0 + n^{-\frac{1}{2}}u)) - \log(\pi(f_0))|\\
&  \leq C n^{-1/2}\sqrt{d} K(d) \leq  C\delta_n.
\end{aligned}
\en
According to Lemma \ref{lem:diffZ}, with high probability, 
\be\label{eq:ratio2}
\Big(\int_{\U} \tilde{Z}_n(u)du\Big)^{-1}\int_{B_K \cap \U}  \Big|Z_n(u) - \tilde{Z}_n(u)\Big| = \mathcal{O}(\delta_n).
\en
This in particular implies that
\be\label{eq:ratio3}
\Big(\int_{\U} \tilde{Z}_n(u)du\Big)^{-1} \int_{B_K\cap \U} Z_n(u)du \leq 1 + o(1).
\en
holds with high probability.
Therefore the theorem follows from \eqref{eq:ratio0}-\eqref{eq:ratio3}.
 
\end{proof} 

\begin{rem}
Using the same arguments as in the proof above, one can prove that any finite moment of the posterior distribution is close to the corresponding moment of the asymptotic normal distribution provided that the dimension $d$ grows much slower than the growth rate described as in Theorem \ref{thm:main}.  Indeed, if $\tilde{\delta}_n := \sqrt{\frac{d}{n}} K(d)^{2+\alpha} \gt 0$ as $n \gt \infty$, then it holds that 
$$
\int_{\F} |f - f_0|^\alpha |\pi_n(f) - \phi(f; \overline{\Delta}_n, n^{-1}\bSig)| df = \mathcal{O}_{p}(\tilde{\delta}_n ).
$$
In particular, by setting $\alpha = 1$, we see that if $\sqrt{\frac{d}{n}} K(d)^{3} \gt 0$, then we obtain the consistency of the posterior mean.

\end{rem}

\section{Application in an Inverse Medium Problem}\label{sec:appl}

\subsection{Forward Model}
Let $\Omega = [0,1]^2 \subset \R^2$. Given two functions $f \in C(\overline{\Omega})$ and $g \in C(\partial \Omega)$, consider the following Dirichlet problem for the elliptic equation on $\Omega$
\be\label{eq:pde}
\bea
& -\Delta u + q u = f & \text{ on } \Omega,\\
& u = g & \text{ on }  \partial \Omega.
\ena
\en
Here we assume that $q\in \A$ where 
$$
\A := \{q\in C(\overline{\Omega}): 0 \leq q_{\min} \leq q(x) \leq q_{\max} < \infty, \quad x\in \Omega \}.
$$
We also assume that $f$ and $g$ are strictly positive on $\Omega$ and $\partial \Omega$ respectively. 
With these assumptions, there exists a unique positive solution $u_q \in C^2(\Omega) \cap C(\overline{\Omega})$ to problem \eqref{eq:pde}. The corresponding inverse problem of interest is the following: Given the solution $u$ on $\Omega$, find the coefficient $q$ on $\Omega$. The problem \eqref{eq:pde} appears as a model problem  for Photo-Acoustic Tomography \cite{BU10} and other multiwave imaging modalities \cite{bal2013hybrid}. The uniqueness of the inverse problem has been proved in \cite{Knowles99} (see also \cite{BU13}).

In practice, the solution can only be measured on a finite set of discrete points. Then the corresponding (finite dimensional) inverse problem is of particular interest to us, that is, to recover the pointwise values of the coefficient at the same observation points. To make this more precise, let us consider the following finite difference approximation to the problem \eqref{eq:pde}:
\be\label{eq:fd}
\begin{cases}
-\frac{u_{i-1, j} + u_{i+1, j} - 4u_{i,j} + u_{i, j-1} + u_{i, j+1}}{h^2} + q_{i, j} u_{i, j}= f_{i, j}, & \forall i, j=1, \cdots, N-1\\
u_{i, 0} = g_{i,0}, u_{i, N} = g_{i,N}, u_{0, j} = g_{0, j}, u_{ N, j} = g_{N, j}, & \forall  i, j=0, \cdots, N.
\end{cases}
\en
Here $f_{i, j} =f(x_i, y_j), q_{i,j} = q(x_i, y_j), g_{i, j} = g(x_i, y_j), i,j=0,\cdots, N$ where $(x_i, y_j) =(i h, j h)$ with the uniform measure size $h = 1/N$. The solution $u_{i,j}$ of the finite difference equations provides an approximation to $u(x_i, y_j)$.  We refer the interested reader to \cite{jovanovic2013analysis} for the convergence analysis of finite difference schemes for PDEs. 

It is more convenient to write the finite difference equation in a matrix form. For doing so, let $\bI$ be a $(N-1)$-dimensional identity matrix and define the $(N-1)^2 \times (N-1)^2$ tridiagonal block matrix $\bA$  and the $(N-1) \times (N-1)$ tridiagonal matrix $\bB$ by
$$
\bA = \left[\begin{array}{cccccc}
\bB & -\bI & & & & \\
-\bI & \bB & -\bI & & &\\
 & & \ddots & \ddots & &\\
 & & & -\bI &\bB & -\bI\\
  & & & & -\bI & \bB\\
\end{array}\right], \quad 
\bB = \left[\begin{array}{cccccc}
4 & -1 & & & & \\
-1 & 4 & -1 & & &\\
 & & \ddots & \ddots & &\\
 & & & -1 &4 & -1\\
  & & & & -1 & 4\\
\end{array}\right].
$$
 We build the solution vector $\bU$, the vector of the right hand side $\bF$ and the vector of the boundary data $\bG$ in the natural row-wise ordering, i.e. 
$$\bea
\mathbf{u} & = \left[u_{1,1}\, \cdots\, u_{N-1,1}\, | \, u_{1,2}\, \cdots\, u_{N-1, 2}\,|\, \cdots \, u_{1, N-1}\, \cdots\, u_{N-1, N-1}\right],\\
\mathbf{f} & = \left[f_{1,1}\, \cdots\, f_{N-1,1}\, |\, f_{1,2}\, \cdots\, f_{N-1, 2}\,| \, \cdots\, f_{1, N-1}\, \cdots\, f_{N-1, N-1}\right],\\
\mathbf{g} & = \Big[g_{0,1} + g_{1,0}\,  g_{0,2}\, \cdots\, g_{0,N-2}\, g_{0,N-1}+g_{1,N}\, |\, g_{2,0}\,\quad \mathbf{0} \quad \, g_{2,N}\, |\, \cdots\\
&\quad \quad \cdots \, |\,  g_{N-2,0}\,\quad \mathbf{0} \quad \, g_{N-2,N}\, |\,
 g_{N-1,0} + g_{N, 1}\,  g_{N,2}\, \cdots\, g_{N,N-2}\, g_{N,N-1}+g_{N-1,N} \Big].
\ena
$$
 Let $\bQ$ be a $(N-1)^2 \times (N-1)^2$ diagonal matrix with diagonal entries made from $\mathbf{q} := \{q_{i, j}\}_{i,j=1}^{N-1}$ in the natural row-wise ordering.  With these notations at hand, the finite difference equation \eqref{eq:fd} can be expressed in a matrix form as
\be\label{eq:ls}
(h^{-2} \bA + \bQ) \mathbf{u} = \mathbf{f} + h^{-2}\mathbf{g}.
\en
Since both $\bA$ and $\bQ$ are positive definite, there exists a unique solution $\mathbf{u}$ to the linear system \eqref{eq:ls}. Hence we can define the forward operator
$$
G : \mathbf{q} \in \F \gt \mathbf{u} \in \R^d \quad \quad \mathbf{u} =  G(\mathbf{q}) = (h^{-2} \bA + \bQ)^{-1} (\mathbf{f} + h^{-2}\mathbf{g})
$$
where $\F := [q_{\min}, q_{\max}]^d$ with $d = (N-1)^2$. Since $\mathbf{u}$ are all viewed as vectors, for convenience we will only use single-index notation $\mathbf{u}_i, i=1,\cdots, d$ instead of $u_{i,j}$ when referencing entries. 

Thanks to the positivity assumption on $f, g$, it follows directly from the discrete maximum principle \cite{Ciarlet70} that there exists  positive constants $C_1\leq C_2$, depending only on $f$ and $g$ such that
\be\label{eq:umm}
C_1 \leq \min_i \mathbf{u}_i \leq |\mathbf{u}|_{\infty}  \leq C_2.
\en
As a consequence, it holds that $|\mathbf{u}|\leq C\sqrt{d}$ for some $C > 0$.
 
Given any fixed $\bq\in \F$, let $\bu_{\bq}$ be the corresponding solution to \eqref{eq:ls}.  Let $\bU_\bq$ be the $d\times d$ diagonal matrix with diagonal entries given by $\bu_\bq$.
The lemma below establishes the differentiability of $G$ and computes explicitly its derivative matrix. 

\begin{lem}\label{lem:dG}
The forward map $\bq \mapsto G(\bq)$ is differentiable at every in $\bq\in\F$ and the derivative matrix $\nabla G(\bq)$ is given by
\be\label{eq:dG}
\nabla G(\bq) = -(h^{-2}\bA + \bQ)^{-1} \bU_{\bq}.
\en

\begin{proof}
Let $\bp$ be any vector in $\R^d$ and let $\eps > 0$. Let $\bu_\eps$ be the solution to \eqref{eq:ls} with $\bQ$ replaced by the diagonal matrix $\bQ_\eps := \text{diag}(\bq + \eps \bp)$, i.e.
\be\label{eq:lseps}
(h^{-2} \bA + \bQ_\eps) \bu_\eps = \mathbf{f}+ h^{-2} \mathbf{g}.
\en
Then subtracting equation \eqref{eq:lseps} by \eqref{eq:ls} leads to
$$
(h^{-2} \bA + \bQ_\eps)(\bu_\eps - \bu_{\bq}) = -\eps (\bQ_\eps - \bQ)\bu_\bq  = -\eps \bU_\bq \bp.
$$
Dividing the above equation by $\eps$ and then letting limit $\eps \lgt 0$, we obtain that 
$$
\bea
\nabla G(\bq) \bp & = \lim_{\eps \lgt 0} \frac{ G(\bq + \eps \bp) - G(\bq)}{\eps}\\
& = \lim_{\eps \lgt 0} \frac{ \bu_\eps - \bu_\bq}{\eps} \\
& = -\lim_{\eps \lgt 0} (h^{-2} \bA + \bQ_\eps)^{-1} \bU_\bq \bp\\
&= -(h^{-2}\bA + \bQ)^{-1} \bU_{\bq}\bp.
\ena
$$
This proves \eqref{eq:dG} since $\bp$ is arbitrary.
\end{proof}
\end{lem}

Notice from \eqref{eq:dG} that the derivative matrix $\nabla G(\bq)$ is symmetric.  The next lemma collects some important properties of $\nabla G$ in the domain $\F$, in correspondence to Assumption \ref{as:G}. 

\begin{lem} \label{lem:dGIP}
\item[(i)] The derivative matrix $\nabla G(\bq)$ is invertible at any $\bq\in\F$.
\item[(ii)] The following holds uniformly with respect to $\bq\in \F$:
\be\label{eq:ii}
d^{-1} \lesssim \lambda_{\min} (\nabla G(\bq)) \leq \lambda_{\max} (\nabla G(\bq)) \lesssim 1.
\en
\item[(iii)] The following estimates hold uniformly with respect to $\bq_1, \bq_2 \in \F$:
\be\label{eq:iii}
|G(\bq_1) - G(\bq_2)|\lesssim |\bq_1 - \bq_2| \lesssim d \times |G(\bq_1) - G(\bq_2)|.
\en
\be\label{eq:iv}
|G(\bq_1) - G(\bq_2) - \nabla G(\bq_2) (\bq_1 - \bq_2)|\lesssim |\bq_1 - \bq_2|^2.
\en

\begin{proof}
\item[(i)] In view of the expression \eqref{eq:dG}, the invertibility of $\nabla G$ follows directly from the positivity of $\bU_\bq$.

\item[(ii)] According to \cite[Example 10.5.1]{Thomas99}, the eigenvalues of the matrix $\bA$ are 
$$
\lambda_{i, j} = 2 \Big(2 - \cos\big(\frac{i\pi}{N}\big) - \cos\big(\frac{j\pi}{N}\big)\Big), \quad i, j = 1, \cdots, N-1.
$$
Hence when $N \gg 1$, or equivalently when $d = (N-1)^2 \gg 1$, we have
\be
\lambda_{\min} (\bA) \sim \mathcal{O}(N^{-2}) \sim \mathcal{O}(d^{-1}) \text{ and } \lambda_{\max}(\bA) = \mathcal{O}(1).
\en
Then it follows from the fact that all the eigenvalues of $\bQ$ are positive and are of order one (since $\bq\in\F$) that
\be\label{eq:eigenAQ}
d^{-1} \lesssim\lambda_{\min} ((h^{-2}\bA + \bQ)^{-1}) \leq \lambda_{\max} ((h^{-2}\bA + \bQ)^{-1}) \lesssim 1.
\en
 Therefore  \eqref{eq:ii} follows from \eqref{eq:dG}, \eqref{eq:eigenAQ} and \eqref{eq:umm}.
 
 \item[(iii)] Let $\bu^1_{\bq}$ and $\bu^2_{\bq}$ be the corresponding solutions to equation \eqref{eq:ls} with $\bq_1$ and $\bq_2$ respectively. Similarly, we can define matrices $\bQ_i$ and $\bU^i_\bq, i=1,2$. Then by using the same arguments as in the proof of Lemma \ref{lem:dG}, one can obtain that 
 \be\label{eq:uq12}
 (h^{-2}\bA + \bQ_1)(\bu^1_{\bq} - \bu^2_{\bq}) = - \bU^2_\bq (\bq_1 - \bq_2).
 \en
which in turn gives 
 $$
 \bq_1 - \bq_2 = -(\bU^2_{\bq})^{-1}  (h^{-2}\bA + \bQ_1)(G(\bq_1) - G(\bq_2)).
 $$
The above equation combining with the bounds in \eqref{eq:eigenAQ} proves \eqref{eq:iii}.
Furthermore, let $\mathbf{v}_{\bq} := \nabla G(\bq_2)(\bq_1 - \bq_2)$. Then according to Lemma \eqref{lem:dG}, we have that
\be\label{eq:vq}
(h^{-2}\bA + \bQ_2) \mathbf{v}_{\bq} = - \bU^2_{\bq} (\bq_1 - \bq_2).
\en
By subtracting equation \eqref{eq:uq12} with equation \eqref{eq:vq}, one sees that
$$
\bea
(h^{-2}\bA + \bQ_1) (\bu^1_{\bq} - \bu^2_{\bq} - \mathbf{v}_\bq) & =  -(\bQ_1 - \bQ_2) \mathbf{v}_\bq \\
& = (\bq_1 - \bq_2)^T \nabla G(\bq_2)(\bq_1 - \bq_2).
\ena
$$
Then \eqref{eq:iv} follows from above equation and eigenvalue estimates \eqref{eq:ii} and \eqref{eq:eigenAQ}.
\end{proof}
\end{lem}

\subsection{Bayesian Inverse Problem} We are interested in estimating the medium parameter $\bq$ on the uniform grid, given the noisy measurements $\bY_n \in \R^d$, where 
$$
\bY_n = G(\bq) + \frac{1}{\sqrt{n}} \eta
$$
with $\mathbf{\eta} \sim N(0,\bI_d)$. To adopt a Bayesian approach to the inverse problem, we put a product prior $\Pi$ on $\bq$ with Lebesgue density 
$$
\Pi(\bq) = \Pi_{i=1}^d \pi_i(\bq_i). 
$$
We assume that the priors on the individual components $\pi_i$ satisfy the following.

\begin{assumptions}\label{as:pii}
\item[(i)] For each $i = 1,2\cdots, d$, $\pi_i$ is supported in $[q_{\min}, q_{\max}]$.

\item[(ii)] There exists $C > 0$ such that for any $q_1, q_2 \in [q_{\min}, q_{\max}]$,
$$
|\log(\pi_i)(q_1) - \log(\pi_i)(q_2)| \leq C
$$

\item[(iii)] Given a $\delta  >0$, there exists $C_\delta > 0$ such that
$$
|\log(\pi_i)(q_1) - \log(\pi_i)(q_2)| \leq C_\delta |q_1 - q_2|
$$
when $|q_1 - q_2| \leq \delta$. 
\end{assumptions}

As we have mentioned in Remark \ref{rem:prior},  if $\pi_i$ satisfies Assumption \ref{as:pii}, then the product prior $\Pi$ constructed as above satisfies Assumption \ref{as:prior}. Let $\Pi_n$ be the resulting posterior, namely
$$
\Pi_n(d\bq) \propto \exp\Big(-\frac{n}{2}|\bY_n - G(\bq)|^2\Big) \Pi(\bq)d\bq.
$$

According to Lemma \ref{lem:dGIP},  the forward map $G$ associated to the medium problem satisfies Assumption \ref{as:G} with $\sigma(d) = d$. Suppose that the data $\bY_n$ is generated from some truth $\bq_0$, which lies strictly inside $\F$. Now we are ready to state the asymptotic normality of the posterior $\Pi_n$ as $n \gt \infty$.
 Let $\bSig = (\nabla G(\bq_0)^T \nabla G(\bq_0))^{-1}$ and let $\Delta_n = \bq_0 + n^{-\frac{1}{2}} \bSig \nabla G(\bq_0)^T \eta$. Then $\Delta_n\sim N(\bq_0, \frac{1}{n}\bSig)$. The following Bernstein-von Mises theorem for the inverse medium problem is simply a restatement of Corollary \ref{cor:bvm-Pi}. 

\begin{thm}\label{thm:bvm-ip}
Let the prior be given by $\Pi(d \bq)  = \Pi_{i=1}^d \pi_i(\bq_i)(d\bq_i)$ with $\pi_i$ satisfying Assumption \ref{as:pii}. If $\delta_n := n^{-\frac{1}{2}} d^3\log d \gt 0$ as $n \gt \infty$, then we have 
$$
d_{TV}(\Pi_n, N(\Delta_n, \bSig_n)) = \mathcal{O}_p(\delta_n)
$$
as $n \gt \infty$.
\end{thm}

\section*{Acknowledgement}
The author thank Professors Richard Nickl, Andrew Stuart and Hendrik Weber for stimulating discussions and suggestions. The author is supported by
EPSRC as part of the MASDOC DTC at the University of Warwick with grant No. EP/HO23364/1.

\bibliographystyle{abbrv}
\bibliography{ref}

\end{document}